\documentclass{article}
\usepackage[utf8]{inputenc}
\usepackage{amsthm, amsfonts, amssymb, latexsym, amsmath}
\usepackage{xcolor}
\usepackage{ulem}

\newcommand{\F}{\mathbb{F}}

\newtheorem{lemma}{Lemma}
\newtheorem{corollary}{Corollary}

\theoremstyle{definition}

\newtheorem{theorem}{Theorem}

\usepackage{refcheck}

\title{Binary sequences derived from differences of consecutive quadratic residues
}
\author{Arne Winterhof$^1$ and Zibi Xiao$^2$\\
~\\
$^1$ Johann Radon Institute for\\
Computational and Applied Mathematics\\
Austrian Academy of Sciences\\
Altenbergerstr.\ 69, 4040 Linz, Austria\\
e-mail: arne.winterhof@oeaw.ac.at\\
$^2$ College of Science\\
Wuhan University of Science and Technology\\ Wuhan 430081, Hubei, China\\
e-mail: xiaozibi@wust.edu.cn
}

\date{}

\begin{document}

\maketitle

\begin{abstract}
    For a prime $p\ge 5$ let $q_0,q_1,\ldots,q_{(p-3)/2}$ be the quadratic residues modulo $p$ in increasing order. We study two $(p-3)/2$-periodic binary sequences $(d_n)$ and $(t_n)$
    defined by $d_n=q_n+q_{n+1}\bmod 2$ and $t_n=1$ if $q_{n+1}=q_n+1$ and $t_n=0$ otherwise,
    $n=0,1,\ldots,(p-5)/2$.
    For both sequences we find some sufficient conditions for attaining the maximal linear complexity $(p-3)/2$.

    Studying the linear complexity of $(d_n)$ was motivated by heuristics of Caragiu et al.
    However, $(d_n)$ is not balanced and we show that a period of $(d_n)$ contains about $1/3$ zeros and $2/3$ ones if $p$ is sufficiently large. In contrast, $(t_n)$ is not only essentially balanced but also all longer patterns of length $s$
    appear essentially equally often in the vector sequence $(t_n,t_{n+1},\ldots,t_{n+s-1})$,
    $n=0,1,\ldots,(p-5)/2$, for any fixed $s$ and sufficiently large $p$.
\end{abstract}



\section{Introduction}

The {\it linear complexity} $L(s_n)$ of a sequence $(s_n)$ over $\F_2$ is the length $L$ of the
shortest linear recurrence
$$s_{n+L}=c_{L-1}s_{n+L-1}+\ldots+c_0s_n,\quad n=0,1,\ldots$$
with coefficients $c_0,\ldots,c_{L-1}\in \F_2$. It is an important measure for the unpredictability and thus suitability of a sequence in cryptography.
For surveys on linear complexity and related measures see \cite{mewi,ni,towi,wi}.

Caragiu et al.\ \cite{catekam} suggested to study the linear complexity of the sequence of the parities of differences of consecutive quadratic residues modulo $p$. In particular,
they calculated the linear complexities for the first $1000$ primes $p\ge 5$.

More precisely,  for a prime $p\ge 5$ we identify the finite field $\F_p$ of $p$ elements with the set of integers $\{0,1,\ldots,p-1\}$. Let $q_0,\ldots,q_{(p-3)/2}$ be the quadratic residues modulo $p$ in increasing order $1=q_0<q_1<\ldots<q_{(p-3)/2}\le p-1$.
We consider the sequence $(d_n)$ of parities of the differences (or sums) of consecutive quadratic residues modulo $p$,
\begin{equation}\label{dndef} d_n=q_n+q_{n+1} \bmod 2,\quad n=0,1,\ldots,(p-5)/2,
\end{equation}
and continue it with period $(p-3)/2$,
\begin{equation}\label{dper} d_{n+(p-3)/2}=d_n,\quad n=0,1,\ldots
\end{equation}
The heuristic of Caragiu et al.\ for the linear complexity of $(d_n)$ shows that among the first $1000$ primes $p\ge 5$ there are $671$ sequences $(d_n)$ with maximal linear complexity $(p-3)/2$.

In Section~\ref{dnlc},
we give some sufficient conditions on $p$ for the maximality of~$L(d_n)=(p-3)/2$.

Balancedness is another desirable feature of a cryptographic sequence, that is, each period should contain about the same numbers of zeros and ones. We show in Section~\ref{imbalance} that the sequence $(d_n)$ contains asymptotically $1/3$ zeros and~$2/3$ ones in each period and is very unbalanced.

Since $(d_n)$ is not balanced, we define a similar $(p-3)/2$-periodic sequence~$(t_n)$ which is essentially balanced and defined by
\begin{equation}\label{tndef} t_n=\left\{\begin{array}{cc} 1, & q_{n+1}=q_n+1,\\
0, & q_{n+1}\ne q_n+1,\end{array}\right.\quad n=0,1,\ldots,(p-5)/2.
\end{equation}

In Section~\ref{pat} we will show that $(t_n)$ is essentially balanced.
Moreover, for fixed length $s$ each pattern $(t_n,t_{n+1},\ldots,t_{n+s-1})=\underline{x}\in \{0,1\}^s$ appears for essentially the same number of $n$ with $0\le n\le (p-3)/2$ provided that $p$ is sufficiently large with respect to $s$.

Finally, we study the linear complexity of $(t_n)$ in Section~\ref{Lt} and provide a sufficient criterion for the maximality of $L(t_n)$. We also prove a lower bound on the
$N$th maximum order complexity of $(t_n)$ which implies a rather moderate but non-trivial and unconditional lower bound on the $N$th linear complexity of~$(t_n)$.

We use the notation $f(n)=O(g(n))$ if $|f(n)|\le cg(n)$ for some absolute constant $c>0$.

\section{Linear complexity of $(d_n)$}\label{dnlc}
Our starting point to determine the linear complexity of a periodic sequence is~\cite[Lemma~8.2.1]{cudire}.
\begin{lemma}\label{cdr}
Let $(s_n)$ be a $T$-periodic sequence over $\F_2$ and
$$S(X)=\sum_{n=0}^{T-1} s_nX^n.$$
Then the linear complexity $L(s_n)$ of $(s_n)$ is
$$L(s_n)=T-\deg(\gcd(X^T-1,S(X))).$$
\end{lemma}
We write the period of the sequence $(d_n)$ in the form
$$T=\frac{p-3}{2}=2^sr$$
with integers $s\ge 0$ and odd $r$. Then we have
$$X^T-1=(X^r-1)^{2^s}.$$
We have to determine $\gcd(X^T-1,D(X))$, where
$$D(X)=\sum_{n=0}^{(p-5)/2} d_nX^n.$$

First we study whether $D(X)$ is divisible by $(X-1)$, that is, we determine the value of $D(1)\in \F_2$.
According to the definition of the sequence $(d_n)$, we get
$$D(1)=\sum_{n=0}^{(p-5)/2} d_n\equiv q_0+2\sum_{i=1}^{(p-5)/2}q_i+q_{(p-3)/2}\equiv 1+q_{(p-3)/2}\bmod 2.$$
Since $-1$ is a quadratic residue modulo $p$ if and only if $p\equiv 1\bmod 4$ and $2$ is a quadratic residue modulo $p$
if and only if $p\equiv \pm 1\bmod 8$,
the largest quadratic residue $q_{(p-3)/2}$ modulo $p$ is
$$q_{(p-3)/2}=\left\{\begin{array}{cc} p-1, & p\equiv 1\bmod4,\\
                                       p-2, & p\equiv 3\bmod 8.
                                       \end{array}\right.$$
In the remaining case $p\equiv 7 \bmod 8$, both $-1$ and $-2$ are
quadratic non-residues. Hence, the largest quadratic residue modulo $p$
is $p-u$ for some $u>2$.
Assume $u=2m$ for some positive integer $m$.
Since $-u$ and $2$ are both quadratic residues modulo $p$, $-m\equiv p-m\bmod p$ is quadratic residue modulo $p$ as well,
a contradiction to the maximality of $p-u$. Hence, $u$ is odd.
So, the largest quadratic residue modulo $p$ is
$$q_{(p-3)/2}=\left\{\begin{array}{ll}
p-1\equiv 0\bmod 2, & p\equiv 1\bmod 4,\\
p-2\equiv 1\bmod 2, & p\equiv 3\bmod 8,\\
p-u\equiv 0\bmod 2, & p\equiv 7\bmod 8.
\end{array}
\right.$$
Thus we have
\begin{equation}\label{S1}
 D(1)=\left\{\begin{array}{ll}
0, & p\equiv 3 \bmod 8,\\
1, & p\not\equiv 3\bmod 8.
\end{array}\right.
\end{equation}


We return now to a general binary sequence $(s_n)$ of period $T$.
The following provides a necessary condition for $S(\beta)=0$ for a primitive~$r$th root of unity $\beta$ in some extension field of $\F_2$.
\begin{lemma}\label{Sbeta}
Let $r$ be an odd prime divisor of $T$ such that $2$ is a primitive root modulo $r$.
Let $\beta$ be any primitive $r$th root of unity in some extension field of~$\F_2$.
If $S(\beta)=0$, then we have
$$\sum_{j=0}^{T/r-1}s_{h+jr}=S(1),\quad h=0,1,\ldots,r-1.$$
\end{lemma}
Proof.
Since $2$ is a primitive root modulo $r$, the cyclotomic polynomial
$$1+X+\ldots+X^{r-1}$$
is irreducible over $\F_2$, and thus the minimal polynomial of $\beta$. In particular we have
$$\beta^{r-1}=\sum_{h=0}^{r-2} \beta^h$$
and $1,\beta,\ldots,\beta^{r-2}$ are linearly independent.
Since $\beta^r=1$ we get
\begin{eqnarray*} S(\beta)&=&\sum_{n=0}^{T-1} s_n\beta^n=\sum_{h=0}^{r-1} \sum_{j=0}^{T/r-1}s_{h+jr}\beta^h\\
&=&\sum_{h=0}^{r-2}
\left(\sum_{j=0}^{T/r-1}s_{h+jr}-\sum_{j=0}^{T/r-1}s_{r-1+jr}\right)\beta^h.
\end{eqnarray*}
Assume $S(\beta)=0$. Then we get
$$\sum_{j=0}^{T/r-1}s_{h+jr}=\sum_{j=0}^{T/r-1}s_{r-1+jr},\quad h=0,1,\ldots,r-2.$$
Hence, since $r$ is odd and
$$S(1)=\sum_{h=0}^{r-1} \sum_{j=0}^{T/r-1}s_{h+jr}=r\sum_{j=0}^{T/r-1}s_{r-1+jr}
=\sum_{j=0}^{T/r-1}s_{r-1+jr},$$
the result follows. ~\hfill $\Box$\\

Now we are ready to prove a sufficient condition on $p$ for $(d_n)$ having maximal linear complexity $L(d_n)=(p-3)/2$.

\begin{theorem}
Let $p=2^{s+1}r+3$ be a prime with $s\in \{0,1\}$ and either $r=1$ or $r$ an odd prime such that $2$ is a primitive root modulo $r$. Then the linear complexity of the sequence $(d_n)$ defined by $(\ref{dndef})$ and $(\ref{dper})$ is maximal,
$$L(d_n)=\frac{p-3}{2}.$$
\end{theorem}
Proof. Since $p=2^{s+1}r+3$ with $s\in \{0,1\}$ and $r$ is odd, we have
$p\not\equiv 3\bmod 8$. It follows from $(\ref{S1})$ that $D(1)=1$.

If $r=1$, that is $T=(p-3)/2=2^s$, we have $X^T-1=(X-1)^{2^s}$,
$\gcd(D(X), X^T-1)=1$ and
$L(d_n)=\frac{p-3}{2}$ by Lemma~\ref{cdr}.

Now let $r$ be an odd prime such that $2$ is a primitive root modulo $r$.
Next we prove that $D(\beta)\neq 0$ for any primitive $r$th root of unity $\beta$.

Assume $D(\beta)=0$.

If $s=0$, we get
$$d_0=d_1=\ldots=d_{r-1}=D(1)=1$$
by Lemma~\ref{Sbeta}.
However, each $n$ with $1\le n\le p-3$ and
$$\left(\left(\frac{n}{p}\right),\left(\frac{n+1}{p}\right),\left(\frac{n+2}{p}\right)\right)
=(1,-1,1),$$
where $\left(\frac{.}{.}\right)$ denotes the Legendre symbol,
corresponds to some $(q_i,q_{i+1})=(n,n+2)$ and thus $d_i\equiv n+n+2\equiv 0\bmod 2$.
By \cite[Proposition~2]{di} there are at least
$$\frac{p}{8}-\frac{\sqrt{p}}{4}-\frac{15}{8}>0,\quad p>25,$$
such $n$, a contradiction for $p>25$. The only remaining primes $p\le 25$ of the form $p=2r+3$ with odd $r>1$ are
$p=13$ and $17$. For $p=13$ we have $q_0=1$ and $q_1=3$, that is, $d_0=0$, a contradiction.
For $p=17$ we get $r=7$ but $2$ is a quadratic residue modulo $7$ and thus not a primitive root modulo $7$.


If $s=1$, we have $p\equiv 7\bmod 8$, $p\ge 23$, and we get from Lemma~\ref{Sbeta}
\begin{equation*} d_0+d_r=d_1+d_{r+1}=\ldots=d_{r-1}+d_{2r-1}=S(1)=1.
\end{equation*}
Hence, $(d_n)$ is balanced.
However, the number of pairs of consecutive quadratic residues is $(p-3)/4$, see for example
\cite[Proposition~4.3.2]{cudire},
and the number of~$n$ with $1\le n\le p-4$ and
\begin{equation}\label{legcond}\left(\left(\frac{n}{p}\right),\left(\frac{n+1}{p}\right),\left(\frac{n+2}{p}\right),
\left(\frac{n+3}{p}\right)\right)=(1,-1,-1,1)
\end{equation}
is at least
$$\frac{p}{16}-\frac{5}{8}\sqrt{p}-\frac{39}{16}>0,\quad p>169,$$
by \cite[Proposition~2]{di}. Hence we have at least
$$\frac{p-3}{4}+\frac{p}{16}-\frac{5}{8}\sqrt{p}-\frac{39}{16}>\frac{p-3}{4}, \quad p>169,$$
different $n$ with $0\le n\le (p-5)/2$ and $d_n=1$, a contradiction for $p>169$.
It remains to check that there is an $n$ satisfying $(\ref{legcond})$ for any prime $p\equiv 7\bmod 8$ for which $r=(p-3)/4$ is a prime and $23\le p<169$, that is,
$p\in\{23,31,47,71,79,127,151,167\}$.
We can delete $p=31,71,127,167$ from this list since for these values of $r=(p-3)/4$ it is easy to verify that $2$ is not a primitive root modulo $r$.
We can choose $n$ from the following table,
$$\begin{array}{c||c|c|c|c}
p & 23 & 47 & 79 & 151\\\hline
n &  9 &  9 &  5 & 5
\end{array}.$$

Thus, we obtain $\gcd(X^T-1,S(X))=1$, and  the result follows.
\hfill $\Box$

\section{Imbalance of $(d_n)$}\label{imbalance}

In this section we show that, for sufficiently large $p$, the sequence $(d_n)$ is imbalanced.
More specifically, about $2/3$ of the sequence elements are equal to~$1$.

\begin{theorem}
Let $N(0)$ and $N(1)$ denote the number of $0$s and $1$s
in a period of the sequence~$(d_n)$, respectively. Then we have
$$N(0)=\frac{p}{6}+O\left(p^{1/2}(\log p)^2\right)$$
and
$$N(1)=\frac{p}{3}+O\left(p^{1/2}(\log p)^2\right).$$
\end{theorem}
Proof.
We first prove a lower bound on $N(1)$.
We need a well known result about the pattern distribution of Legendre symbols.

For $s\ge 1$ and $\varepsilon_1, \varepsilon_2, \cdots, \varepsilon_s\in \{-1,1\}$, set
$$ N(\varepsilon_1, \cdots, \varepsilon_s)
=\left| \left\{j=1,2,\ldots,p-s: \left(\frac{j+i}{p}\right)=\varepsilon_{i+1},~i=0,\ldots, s-1\right\}\right|.
$$
From \cite[Proposition~2]{di} we get for $s\geq 3$,
\begin{equation}\label{e2}
N(\varepsilon_1, \cdots, \varepsilon_s)=\frac{p}{2^{s}}+O\left(sp^{1/2}\right).
\end{equation}
Note that $(\ref{e2})$ is also true for $s=1$, since
we have each $(p-1)/2$ quadratic residues and non-residues modulo $p$, and for $s=2$, see for example
\cite[Proposition~4.3.2]{cudire}.

For a non-negative integer $k$, let $N_k$ denote the number of $j$ with $1\le j\le p-2-k$
satisfying
\begin{equation}\label{jk}\left(\left(\frac{j}{p}\right),\left(\frac{j+1}{p}\right),\cdots,\left(\frac{j+k}{p}\right),\left(\frac{j+k+1}{p}\right)
\right)=(1,-1,\cdots,-1,1).
\end{equation}
Each pair $(j,k)$ satisfying $(\ref{jk})$ corresponds to an $n$ with $(q_n,q_{n+1})=(j,j+k+1)$, that is,
$d_n\equiv q_n+q_{n+1}\equiv k+1\bmod 2$.
Hence for any positive integer $m$,
$$N(1)\ge \sum_{k=0}^m N_{2k}=\frac{p}{4}\sum_{k=0}^m4^{-k}+O\left(m^2p^{1/2}\right)$$
and
$$N(0)\ge \sum_{k=0}^m N_{2k+1}=\frac{p}{8}\sum_{k=0}^m4^{-k}+O\left(m^2p^{1/2}\right)$$
by $(\ref{e2})$.
Choosing $m=\lfloor \log p\rfloor$
we get
\begin{eqnarray*}
 N(1) &\ge& 
   \frac{p}{3}\left(1-\left(\frac{1}{4}\right)^{m+1}\right)+O\left(m^2p^{1/2}\right)\\
    &=& \frac{p}{3}+O\left(p^{1/2}(\log p)^2\right)
\end{eqnarray*}
and
$$N(0)\ge \frac{p}{8}\sum_{k=1}^m 4^{-k}+O(m^2p^{1/2})
=\frac{p}{6}+O\left(p^{1/2}(\log p)^2\right).$$
Now since $N(0)+N(1)=(p-3)/2$ we get
$$N(0)= \frac{p}{6}+O(p^{1/2}(\log p)^2) \quad \mbox{and}\quad N(1)= \frac{p}{3}+O(p^{1/2}(\log p)^2).$$
Therefore, the sequence $(d_n)$ is imbalanced for sufficiently large $p$.~\hfill $\Box$

\section{Pattern distribution of $(t_n)$}\label{pat}

The number $N(1)$ of $1$s in a period of the sequence $(t_n)$ defined by $(\ref{tndef})$ is equal to
the number of elements of the set
$$ \left\{j=1,2,\ldots,p-2:  \left(\frac{j}{p}\right)=\left(\frac{j+1}{p}\right)=1\right\}.$$
Then it follows from \cite[Proposition 4.3.2]{cudire} that
\begin{equation}\label{e5}
N(1)=\left\{\begin{array}{ll}
(p-3)/4, & p\equiv 3\,(\bmod\,4),\\
(p-5)/4, & p\equiv 1\,(\bmod\,4).
\end{array}\right.
\end{equation}
So this sequence is balanced when $p\equiv 3\bmod 4$ and almost balanced when $p\equiv 1\bmod 4$.

Now we consider longer patterns.
\begin{theorem}\label{thmpatt}
For a prime $p\ge 5$ let $(t_n)$ be the $(p-3)/2$-periodic sequence defined by $(\ref{tndef})$.
For any positive integer $s$ and any pattern $\underline{x}=(x_0,\ldots,x_{s-1})\in \{0,1\}^s$
the number $N_s(\underline{x})$ of $n$ with $0\le n\le (p-5)/2$ and
$$(t_n,t_{n+1},\ldots,t_{n+s-1})=\underline{x}$$
satisfies
$$N_s(\underline{x})=\frac{p}{2^{s+1}}+O\left(sp^{1/2}(\log p)^{s+1}\right).$$
\end{theorem}
Proof. Each pattern of Legendre symbols
\begin{eqnarray*}
&&\left(\left(\frac{j}{p}\right),\left(\frac{j+1}{p}\right),\ldots,\left(\frac{j+k_0+\ldots+k_{s-1}+s}{p}\right)\right)\\
&=&(1,\underbrace{-1,\ldots,-1}_{k_0},1,\underbrace{-1,\ldots,-1}_{k_1},1,\ldots,
1,\underbrace{-1,\ldots,-1}_{k_{s- 1}}, 1),
\end{eqnarray*}
$j=1,2,\ldots,p-1-k_0-\ldots-k_{s-1}-s$, corresponds to a pattern
$(t_n,\ldots,t_{n+s-1})$
with
$$t_{n+i}=\left\{\begin{array}{cc}1, & k_i=0,\\ 0, & k_i>0,\end{array}\right.\quad i=0,\ldots,s-1,$$
for some $n$ with $0\le n\le \frac{p-5}{2}-s$.
Assume $$m\ge \max\{1,k_0,k_1,\ldots,k_{s-1}\}.$$
Then the number of such $j$ is
$$\frac{p}{2^{s+1+k_0+\ldots+k_{s-1}}}+O(smp^{1/2})$$
by $(\ref{e2})$.

Assume that the pattern $\underline{x}$
contains $r$ zeros.
Then for $r\ge 1$ we have
\begin{eqnarray*}
N_s(\underline{x})&\ge& \frac{p}{2^{s+1}}\sum_{\ell_1,\ldots,\ell_r=1}^m 2^{-(\ell_1+\ldots+\ell_r)}
+O\left(sm^{r+1}p^{1/2}\right)\\
&=&\frac{p}{2^{s+1}}(1-2^{-m})^r+ O\left(sm^{r+1}p^{1/2}\right).
\end{eqnarray*}
Choosing $m=\lfloor \log p\rfloor-1$
we get
\begin{equation}\label{Nsabs}
N_s(\underline{x})\ge \frac{p}{2^{s+1}}+O\left(sp^{1/2}(-1+\log p)^{r+1}\right).
\end{equation}
For $r=0$ we get
\begin{equation*}
N_s(\underbrace{1,1,\ldots,1}_ s)=\frac{p}{2^{s+1}}+O\left(sp^{1/2}\right).
\end{equation*}
Using
\begin{eqnarray*} N_s(\underline{x})&\le&\frac{p-3}{2}-\sum_{\underline{y}\in \F_2^s\setminus\{\underline{x}\}}N_s(\underline{y})\\
&\le& \frac{p-3}{2}-(2^s-1)\frac{p}{2^{s+1}}+O\left(sp^{1/2}\sum_{r=0}^{s}{s\choose r}(-1+\log p)^{r+1}\right)\\
&=&\frac{p}{2^{s+1}}+O\left(sp^{1/2}(\log p)^{s+1}\right)
\end{eqnarray*}
we get the result.
 \hfill $\Box$\\

Using \cite[Theorem~3]{masa} instead of \cite[Proposition~2]{di} we get a local analog of Theorem~\ref{thmpatt} exactly the same way.
\begin{corollary}
For a prime $p\ge 5$ let $(t_n)$ be the $(p-3)/2$-periodic sequence defined by $(\ref{tndef})$.
For any positive integer $s$, any $N$ with $1\le N\le (p-5)/2$ and any pattern $\underline{x}=(x_0,\ldots,x_{s-1})\in \{0,1\}^s$
the number $N_s(\underline{x},N)$ of $n$ with $0\le n\le N-1$ and
$$(t_n,t_{n+1},\ldots,t_{n+s-1})=\underline{x}$$
satisfies
$$N_s(\underline{x},N)=\frac{N}{2^{s+1}}+O\left(sp^{1/2}(\log p)^{s+2}\right).$$
\end{corollary}
We also get an analog of the lower bound $(\ref{Nsabs})$,
\begin{equation}\label{new}
N_s(\underline{x},N)\ge \frac{N}{2^{s+1}}+O\left(sp^{1/2}(\log p)^{r+2}\right),
\end{equation}
where $r$ is the number of zeros of $\underline{x}\in \{0,1\}^s$.

\section{Linear complexity of $(t_n)$}\label{Lt}
In this subsection we discuss the linear complexity of the sequence $(t_n)$. We now put $$T(X)=\sum_{n=0}^{(p-5)/2} t_nX^n.$$
According to $(\ref{e5})$, the number $N(1)$ of $1s$ in a period of $(t_n)$
is equal to $(p-3)/4$ if $p\equiv 3\bmod 4$ and $(p-5)/4$ if $p\equiv 1\bmod 4$. Thus,
$$T(1)=\sum_{n=0}^{(p-5)/2}t_n
= \left\{\begin{array}{ll} 1, & p\equiv \pm 1\bmod 8,\\ 0,& p\equiv \pm 3\bmod 8.\end{array}\right.$$

For the case $p\equiv 1 \bmod 4$, the period $T=r=(p-3)/2$ of the
sequence $(t_n)$ is an odd number. If we suppose that $r$ is
a prime such that $2$ is a primitive root modulo $r$, then
Lemma~\ref{Sbeta} implies either $T(\beta)\ne 0$ or $t_h=T(1)$ for all $h$.
Now $2$ can be only a primitive root modulo $r$ if it is not a square modulo $r$, that is,
$r\equiv \pm 3\bmod 8$ and thus $p\equiv 9,13\bmod 16$, in particular, we have $p\ge 13$
and~$(t_h)$ is not constant by $(\ref{e5})$. Hence, $T(\beta)\ne 0$ for any primitive $r$th root of unity~$\beta$.
We obtain the following result.

\begin{theorem}
Let $p$ be a prime with  $p\equiv 9$ or $13 \bmod 16$ such that $r=\frac{p-3}{2}$ is an odd prime and $2$ is a primitive root modulo $r$.
Then the linear complexity~$L(t_n)$ of the sequence~$(t_n)$
defined by $(\ref{tndef})$ is
$$L(t_n)=\left\{\begin{array}{cc} \frac{p-3}{2}, & p\equiv 9 \bmod 16,\\
\frac{p-5}{2}, & p\equiv 13 \bmod 16.\end{array}\right.$$
\end{theorem}

The {\it maximum order complexity}  $M(s_n)$ of a binary sequence $(s_n)$ is the smallest positive integer~$M$ with
$$ s_{n+M}=f(s_{n+M-1},\ldots,s_n),\quad n=0,1,\ldots,
$$
 for some mapping $f:\F_2^M \mapsto \F_2$.
Obviously, we have
 $$L(s_n)\ge M(s_n)$$
 and each lower bound on $M(t_n)$ is also a lower bound on $L(t_n)$.
In particular we have the trivial lower bound
 $$L(t_n)\ge M(t_n)\ge \frac{\log((p-3)/2)}{\log 2},$$
see \cite[Proposition~3.2]{ja}.

For a positive integer $N$ the {\it $N$th maximum order complexity} $M(s_n,N)$ is the local analog of $M(s_n)$, that is, the smallest $M$ with
$$ s_{n+M}=f(s_{n+M-1},\ldots,s_n),\quad n=0,1,\ldots,N-M-1,
$$
for some $f$.
We prove also a lower bound on $M(t_n,N)$ which is nontrivial for~$N$ of order of magnitude at least $p^{1/2}\log^4 p$.

\begin{theorem}\label{maxord}
For the $N$th maximum order complexity $M(t_n,N)$ of the sequence~$(t_n)$
defined by
$(\ref{tndef})$ we have
$$M(t_n,N)\ge \frac{\log (N/p^{1/2})}{\log 2}- \frac{4\log\log p}{\log 2}+O(1),\quad N=1,2,\ldots,(p-5)/2.$$
\end{theorem}
Proof. For $s\ge 1$ and $x\in\{0,1\}$ the number $G_{s,x}(N)$
of~$n$ with $0\le n\le N-s$ satisfying
\begin{equation}\label{pm}
(t_n,t_{n+1},\ldots,t_{n+s-2},t_{n+s-1})=(\underbrace{1,1,\ldots,1}_{s-1},x )
\end{equation}
satisfies
$$G_{s,x}(N)\ge\frac{N}{2^{s+1}}+O\left(sp^{1/2}(\log p)^3\right)=\frac{N}{2^{s+1}}+O\left(p^{1/2}(\log p)^4\right)\quad \mbox{for }s\le \log p$$
by $(\ref{new})$.
Hence, there is a constant $c>0$ such that for
$$s\le  \frac{\log (N/p^{1/2})}{\log 2}-\frac{4\log\log p}{\log 2}-c$$
we have $G_{s,x}>0$ for $x\in\{0,1\}$ and both patterns in $(\ref{pm})$ of length $s$ appear at least once.
Assume
\begin{equation}\label{M} M\le \frac{\log (N/p^{1/2})}{\log 2}-\frac{4\log\log p}{\log 2}-c
\end{equation}
and that there is a recurrence of the form
\begin{equation}\label{trek} t_{n+M}=f(t_{n+M-1},\ldots,t_n),\quad n=0,1,\ldots,N-M-1.
\end{equation}
However, there are $n_1$ and $n_2$ with $0\le n_1< n_2\le N-1-M$
and
$$t_{n_1+i}=t_{n_2+i}=1,\quad i=0,\ldots,M-1,\quad t_{n_1+M}\ne t_{n_2+M},$$
a contradiction  to $(\ref{trek})$. Hence,  $(\ref{M})$ is not true and the result follows.\hfill $\Box$\\

Remark.
Theorem~\ref{maxord} is in good correspondence to the result of \cite{ja} that the maximum order complexity of a random sequence of length $N$ is of order of magnitude $\log N$.

For some recent papers on the maximum order complexity see \cite{ge,iwi,lu,pe,suwi0,suwi,szlh,xzxlj}.

The {\it correlation measure} $C_2(s_n)$ {\it of order $2$} of a sequence $(s_n)$ of length $N$
is defined by
$$C_2(s_n)=\max_{M,d_1,d_2}\left|\sum_{n=0}^{M-1}(-1)^{s_{n+d_1}+s_{n+d_2}}\right|,$$
where the maximum is taken over all integers $M,d_1,d_2$ with $0\le d_1< d_2\le N-M$.
There exist $d_1$ and $d_2$ with $0\le d_1<d_2$ with
$s_{n+d_1}=s_{n+d_2}$ for $n=0,1,\ldots,M(s_n)-2$
and we get
\begin{equation}\label{cm} C_2(s_n)\ge M(s_n)-1.
\end{equation}
A large correlation measure $C_2(s_n)$ of order $2$ is undesirable for cryptographic applications since
the expected value of $C_2(s_n)$ is of order of magnitude
$$N^{1/2}(\log N)^{1/2},$$
see~\cite{al}, and a cryptographic sequence should not be distinguishable from a random sequence.

These results on expected values and $(\ref{cm})$ suggest that a good cryptographic sequence of length $N$ should have maximum order complexity of order of magnitude between~$\log N$ and $N^{1/2+\varepsilon}$.

\section{Conclusion}
We showed that the sequence $(d_n)$ of the parities of differences of quadratic residues modulo $p$ is very unbalanced. Hence, $(d_n)$ is, despite of a high linear complexity (at least in some cases), not suitable in cryptography. We introduced an alternative sequence $(t_n)$ which is not only balanced but also longer patterns appear essentially equally often. Moreover, we proved that $(t_n)$ has in some cases a very high linear complexity and obtained a moderate but nontrivial lower bound on the $N$th maximum order complexity of $(t_n)$. All these results indicate that $(t_n)$ is an attractive candidate for applications in cryptography.

\section*{Acknowledgments}
The first author is partially supported by the Austrian Science Fund FWF Project P 30405-N32.
The second author is supported by the Chinese Scholarship Council.

We wish to thank the anonymous referees for their careful study of our paper and their very useful comments.

\end{document}